\documentclass[12pt, a4paper, titlepage]{article}
\usepackage{amsthm,amsmath, amsfonts, amssymb}
\begin{document}

\begin{center}
\huge \bf
On the Chekhov-Fock coordinates of dessins d'enfants.\\
\Large \bf
G. Shabat and V.Zolotarskaia
\end{center}
\begin{center}
\large \bf
Introduction\\
\end{center}
\large
There are several ways to associate a complex structure to a ribbon graph.
The examples are provided by the construction of Kontzevich \cite{Kon}, the construction of
Penner \cite{Pen}, the construction of dessins d'enfants \cite{Chp}. In the
first two constructions the ribbon graph is considered with the additional
structure - a number is associated to each edge. Varying these
parametres we obtain different riemann surfaces. The statement is that this
way we cover a cell of the corresponding moduli space. In the construction of
dessins d'enfants a single riemann surface is associated to each graph.
(Usually we say \it dessin d'enfants, \rm which in this paper means exactly the same as \it ribbon graph).
\rm We call it \it the Grothendieck model of a ribbon graph.

\rm
The
problem is: which parametres for the edges of graph in the first two
constructions should be chosen
to obtain its Grothendiek model? This problem
was solved in \cite{K-CP}, \cite{P-CP} and it turns out that these
parameters should all be equal to $1$.

The goal of the present paper is to discuss one more such construction -
that of Chekhov-Fock and to prove that putting all
parametres of the edges of the graph
in this construction equal
to $0$, we obtain the Grothendieck model of this graph.
\begin{center}
\bf Cartography
\end{center}
Let $\Gamma$ be any trivalent ribbon graph, that is a graph with the valencies
of all the vertices equal to $3$ and with the given cyclic order
on the origins of the edges in each vertice.
Let $E$ be the set of the
oriented edges of $\Gamma$. We have the so-called \it cartography group \rm $C_2^+$ acting on $E$.
This group is generated by the elements $\rho_0$ and $\rho_1$. The element $\rho_0$
turnes the current edge contraryclockwise around the origin
of the edge, $\rho_1$ changes the orientation of the edge (see \cite{Chp}).
Formally we can write
$$ \it C_2^+:=<\rho_0, \rho_1 | \rho_1^2=1>.$$
Having in mind the trivalency of our graph we define $$ \it C_2^+[3]=<\rho_0, \rho_1 |
\rho_1^2=\rho_0^3=1>$$ also acting on $E$.

Fix $\epsilon \in E$.  Let $\it B(E,\epsilon)$ be the borel subgroup of $\it
C_2^+[3]$ corresponding to the edge $\epsilon$ (see \cite{Chp}), that is $$\it B(E,\epsilon)=\{w \in
C_2^+[3]|w\epsilon=\epsilon\}.$$
\begin{center}
\bf The Chekhov-Fock construction
\end{center}

We also have the additional structure:
$$z:E \longrightarrow \bf \mathbb{R}, \it \quad z(\rho_1\gamma) =z(\gamma) \quad \forall \gamma \in \Gamma.$$
Now given $\Gamma, \it z, \epsilon$ we define the map
$$CHF: \it C_2^+[3] \longrightarrow PSL_2(\mathbb {R})$$
inductively, setting
$$CHF(1)=1,$$
$$CHF(\rho_0 w)=L \times CHF(w)$$
$$CHF(\rho_1 w)=X_{z(w\epsilon)} \times CHF(w)$$
where
\large
\\
$X_a=
 \left(\begin{array}{cc}
           0& -e^{a/2} \\
           e^{-a/2}& 0
   \end{array}\right)
$
\\
$L=
 \left(\begin{array}{cc}
           0 & 1\\
           -1 & -1\\
   \end{array}\right)
$
\bf
\\
Notation. \it $chf:=CHF |_{B(E,\epsilon)}$
\\
\bf Statement 1. \it The map $chf$ is a homomorphism.

\rm This statement, as most of the other statements of present paper, becomes obvious after thinking about it for a
while, but here is the formal
\\
\bf Proof:
\rm We should prove that  $$\it chf(w_2w_1)= chf(w_2)
chf (w_1)$$
\rm for any  $\it w_1, w_2\in {B(E,\epsilon)}$.
\\
It is sufficient to show that $$\it CHF(w_2w_1)=
CHF(w_2)CHF(w_1)$$
for any $\it w_2\in C_2^+[3], w_1\in
{B(E,\epsilon)}$.
\\
We will do it using induction over length of $\it w_2$.
\\
If the length of
$\it w_2$ is $1$ then either $\it w_2=\rho_0$,
\\
or $\it w_2=\rho_1$.
\\
In the first case we have $$\it CHF(w_2w_1)=$$ $$\it CHF(\rho_0w_1)= \qquad
( by\; the\; definition \;of \;CHF)$$ $$\it L\times CHF(w_1)=  \qquad ( by  \;
the  \; definition \;  of \;  CHF)$$ $$\it CHF(\rho_0)CHF(w_1)=$$
$$\it CHF(w_2)CHF(w_1)$$
\\
In the second case we have $$\it CHF(w_2w_1)=$$ $$\it CHF(\rho_1w_1)=
\qquad (by  \; the \; definition \; of  \; CHF)$$
$$\it X_{z(w_1\epsilon)}\times CHF(w_1)=
\qquad ( by \; the \; definition \;of \;B(E,\epsilon) )$$
$$\it X_{z(\epsilon)}\times CHF(w_1)= \qquad (by\; the\; definition\; of\; CHF)$$
$$\it CHF(\rho_1)CHF(w_1)=$$
$$\it CHF(w_2)CHF(w_1)$$
\\
In the general case if $\it w_2=\rho_0w_3$ we have
$$\it CHF(w_2w_1)=$$
$$\it CHF(\rho_0w_3w_1)= \qquad (by\; the\; definition\; of\; CHF)$$
$$\it L\times CHF(w_3w_1)= \qquad (by\; the\; induction)$$
$$\it L\times CHF(w_3)CHF(w_1)=\qquad (by\; the\; definition\; of\; CHF)$$
$$\it CHF(\rho_0w_3)CHF(w_1)=$$
$$\it CHF(w_2)CHF(w_1)$$
\\
Finaly if $\it w_2=\rho_1w_3$ we have
$$\it CHF(w_2w_1)=$$
$$\it CHF(\rho_1w_3w_1)=\qquad (by\; the\; definition\; of\; CHF)$$
$$\it X_{z(w_3w_1\epsilon)}\times CHF(w_3w_1)=\qquad (by\;
the\;
definition\; of\; B(E,\epsilon))$$
$$\it X_{z(w_3\epsilon)}\times CHF(w_3)CHF(w_1)=\qquad (by\; the\;
definition\; of\; CHF)$$
$$\it CHF(\rho_1w_3)CHF(w_1)=$$
$$\it CHF(w_2)CHF(w_1)$$
\rm
\\which finishes the proof.
$\blacksquare$
\begin{center}
\bf The Chekhov-Fock net
\end{center}

Let us denote the image of $chf$ by $\Delta(\Gamma, z, \epsilon)$.
To explain the properties of $\Delta(\Gamma, z, \epsilon)$ we need the
notion of \it Chekhov-Fock net. \rm It is some ideal triangulation of the upper
plane $\mathcal {H}$ with the numbers associated to the edges of its dual trivalent graph.

 Here is the construction of Chekhov-Fock net, associated to the graph $\Gamma$. The first triangle of
Chekhov-Fock net will be the ideal triangle $T_0$ with the vertices in
$-1,0$ and $\infty$. Denote $[x, y]$ the Lobachevsky line joining $x$ and $y$.
Let us call $\epsilon^*$ the edge of the trivalent graph, intersecting the edge
$[0,\infty]$ of the triangle (we need the orientation of the edges, so let us
say that the origin of $\epsilon^*$ which is situated inside $T_0$ is its "beginning".
The number, associated to this edge of the graph, will be $\it z^*(\epsilon^*)
=z(\epsilon)$. Now (using the fact that $C_2^+[3]$ acts on our trivalent
graph and that this action is transitive) all the numbers associated to all the edges of the graph
can be determined in the following way: $\it
z^*(w\epsilon^*)=z(w\epsilon)$, where $\it w \in C_2^+[3]$.
Now we have to explain, how to obtain inductively all the other triangles. For example
we have constructed the triangle $T$ with the edges $a,b$ and $c$ and we want to
construct the triangle, which will also contain edge $c$. Let $\alpha \in PSL_2(\mathbb {R})$ be
a transformation of $\mathcal{H}$, which takes triangle $T$ to the triangle
$T_0$, so that $\alpha(c)=[0, \infty]$. Then the
$(\alpha^{-1}X_{z^*(c)}\alpha)T$ is the desired triangle.
\\
\bf Fact. $\it \forall w \in C_2^+[3] \quad z^*((CHF(w))^{-1}\epsilon^*)=z(w\epsilon)$
\\
\bf Statement 2. \it For any $\it w \in C_2^+[3]$ $\it \exists$ triangle $T$ of Fock-Chehov net
  such that $(CHF(w))^{-1}T_0=T$.
\\
\bf Proof: \rm Let us use the induction over the length of $\it w$. if this
length is $0$ then the statement is trivial. Now if $\it w=\rho_0w_1$,
$$\it (CHF(w))^{-1}T_0= \qquad (by\; the\; definition\; of\; CHF)$$
$$(L \times CHF(w_1))^{-1}T_0=$$
$$(CHF(w_1))^{-1}L^{-1}T_0=$$
$$(CHF(w_1))^{-1}T_0=T.$$
Now consider the case $\it w=\rho_1w_1$, let $(CHF(w))^{-1}T_0=T$.
$$ \it (CHF(w))^{-1}T_0= \qquad (by\; the\; definition\; of\; CHF)$$
$$(X_{z(w\epsilon)} \times CHF(w_1))^{-1}T_0=$$
$$(CHF(w_1))^{-1}(X_{z(w\epsilon)})^{-1}T_0=$$
$$((CHF(w_1))^{-1}(X_{z(w\epsilon)})^{-1}CHF(w_1))(CHF(w_1))^{-1}T_0=$$
$$((CHF(w_1))^{-1}(X_{z(w\epsilon)})^{-1}CHF(w_1))T= \qquad (see \; \bf Note \rm \; )$$
$$((CHF(w_1))^{-1}(X_{z^*((CHF(w))^{-1}\epsilon^*)})^{-1}CHF(w_1))T$$
And the last expression by the definition gives us the
next triangle of Chekhov-Fock net.
$\blacksquare$
\bf
\\
Corollary: \it $\Delta ( \Gamma, z, \epsilon)$ is a fuchsian group.

\rm Triangle $T_0$ is its fundamental domain.
\begin{center}
\bf The Chekhov-Fock coordinates of dessins d'enfants
\end{center}

\rm
Now let us put $z\equiv 0$.
\bf
\\
Statement 3: The riemann surface \it $\mathcal{H} \diagup \Delta ( \Gamma, 0, \epsilon)$ is equivalent to the Grothendieck model of $\Gamma$.
\bf
\\
Proof: \rm
If $z\equiv 0$ then $\Delta(\Gamma, z, \epsilon)=\Delta ( \Gamma, 0, \epsilon)$ is the subgroup of
$PSL_2[\mathbb{Z}]$
Let us first consider the case of normal $B(E,\epsilon)$, (we call such
dessin regular).
 Consequently $\Delta(\Gamma, 0, \epsilon)$ is normal in $PSL_2[Z]$.
So we can factorize $\mathcal{H} \diagup \Delta (\Gamma, 0, \epsilon )$ by
$PSL_2[\mathbb{Z}]$. Now if identify $\mathcal{H}/PSL_2[\mathbb{Z}]$ and Riemann sphere, we get the
function $\beta_\Gamma$
defined on $\mathcal{H} / \Delta (\Gamma, 0, \epsilon )$ with the only critical
values $0,1$ and $\infty$, that is Belyi function corresponding to this dessin.

Now if we take any graph $\Gamma$ and if $D$ is the corresponding
dessin d'enfant, there exists regular dessin d'enfant
$D'$  and the covering $\chi: D' \longrightarrow D$. Let us denote the graph
corresponding to this dessin by $\Gamma'$. Now we can define
function $\beta_\Gamma$ defined on $\Sigma / \Delta (\Gamma, 0, \epsilon )$ such that
$\beta_\Gamma\circ\chi=\beta_\Gamma'$
\begin{center}
\bf
Examples.\
\end{center}
\bf
Notation
\rm
By the case $<a_1, \dots ,a_n\mid b_1, \dots , b_m>$ we denote graph with the
valences of vertices equal to $a_1, \dots , a_n$ and the valences of
vertices of dual graph equal to $b_1, \dots , b_m$.
\bf
\\
Example 1. Case$\langle3,3|2,2,2\rangle$
\rm
\large
Let the numbers, corresponding to the edges $A,B,C$, be $a,b,c$,
the chosen oriented edge be B with the orientation, so that edge $A$ is to
the left of the end of $B$ and edge $C$ - to the right.
Then the generators of the corresponding fuchsian group are
\\
$\gamma_1=X_bRX_cR$
\\
$\gamma_2=RX_aRX_b$
\\
$\gamma_3=LX_cRX_aL$
\\
with the relation $\gamma_3\gamma_2\gamma_1=1$.

Using the fact that $\gamma_i$ are parabolic we get
$
\left\{\begin{array}{lcl}
a+b&=&0\\
a+c&=&0\\
b+c&=&0\\
\end{array}\right.
$

This system has the only solution
$
\left\{\begin{array}{lcl}
a&=&0\\
b&=&0\\
c&=&0\\
\end{array}\right.
$

So we get the subgroup of $PSL_2(\bf R)$ with generators:

$\gamma_1= \left(\begin{array}{cc}
           1& 0 \\
           2& 1
   \end{array}\right)$
, the invariant point is $0$.

$\gamma_2= \left(\begin{array}{cc}
           1& -2 \\
           0& 1
   \end{array}\right)$
, the invariant point is $\infty$.

$\gamma_3= \left(\begin{array}{cc}
           -1& -2 \\
           2& 3
   \end{array}\right)$
, the invariant point is $-1$.
\\
with the relation $\gamma_3\gamma_2\gamma_1=1$.
\bf

Example 2. Case$\langle3,3,3,3|3,3,3,3\rangle - tetrahedron $.

\rm
Let the edges $A,B,C$ follow clockwise around some vertice, edge $D$ be the
opposite to the edge $A$, $E$ to $B$, $F$ to $C$.
Taking $a,b,c,d,e,f$ by numbers, corresponding to the edges $A,B,C,D,E,F$
and using the fact that all elements of fuchsian group are parabolic we get
the system
$
\left\{\begin{array}{lcl}
a+e+c&=&0\\
b+d+c&=&0\\
a+f+b&=&0\\
e+f+d&=&0\\
\end{array}\right.
$

The solution of this system is
$
\left\{\begin{array}{lcl}
a&=&d\\
b&=&e\\
c&=&f\\
a+b+c&=&0\\
\end{array}\right.
$

So we get the family of riemann surfaces with two parameters. The riemann
surface coresponding to the dessin can be obtained putting
$a=b=c=0$, because of its symmetry.

So we get the following matrices
\\
$\gamma_1= \left(\begin{array}{cc}
           1& -3 \\
           0& 1
   \end{array}\right)$
, the invariant point is $\infty$
\\
$\gamma_2=\left(\begin{array}{cc}
           2& 3 \\
           -3& -4
   \end{array}\right) $
, the invariant point is $-1$
\\
$\gamma_3=\left(\begin{array}{cc}
           1& 0 \\
           3& 1
   \end{array}\right) $
, the invariant point is $0$
\\
$\gamma_4= \left(\begin{array}{cc}
           4& -3 \\
            3& -2
   \end{array}\right)$
, the invariant point is $1$
\\
with the relation $\gamma_4\gamma_3\gamma_2\gamma_1=1$.

\bf Example 3. Case$\langle3,3,3,3,3,3,3,3|4,4,4,4,4,4\rangle$
\rm
Using the usual technique of Fock-Chehov we get the generators of the
fuchsian group
\\
$\gamma_1= \left(\begin{array}{cc}
           1& 0 \\
           4& 1
   \end{array}\right)$
, the invariant point is $0$
\\
$\gamma_2= \left(\begin{array}{cc}
           5& -4 \\
           4& -3
   \end{array}\right)$
, the invariant point is $1$
\\
$\gamma_3=\left(\begin{array}{cc}
           1& -4 \\
           0& 1
   \end{array}\right) $
, the invariant point is $\infty$
\\
$\gamma_4=\left(\begin{array}{cc}
           7& 16 \\
           -4& -9
   \end{array}\right) $
, the invariant point is $-2$
\\
$\gamma_5=\left(\begin{array}{cc}
           3& 4 \\
           -4& -5
   \end{array}\right) $
, the invariant point is $-1$
\\
$\gamma_6= \left(\begin{array}{cc}
           7& 4 \\
           -16& -9
   \end{array}\right)$
, the invariant point is $-1/2$
\\
with the relation $\gamma_6\gamma_5\gamma_4\gamma_3\gamma_2\gamma_1=1$.
\\
Also we can find the automorfism of this picture:
$\alpha= \left(\begin{array}{cc}
           1& 1 \\
           0& 1
   \end{array}\right)$

Factorizing this child's picture by this automorfism we get the
another picture:
\\
\bf Example 4. Case$\langle3,3|4,1,1\rangle$
\rm
with generators
\\
$\gamma_1= \left(\begin{array}{cc}
           1& 1 \\
           0& 1
   \end{array}\right)$
, the invariant point is $\infty$
\\
$\gamma_2= \left(\begin{array}{cc}
           1& 0 \\
           4& 1
   \end{array}\right)$
, the invariant point is $0$


\begin{thebibliography}{10}

\bibitem{FC}
\large
 \rm V. V. Fock and L. O. Chekhov, \bf
Quantum Mapping Class Group, Pentagon Relation, and Geodesics,
 \it Proceedings of the Steklov Institute of Mathematics, Vol.226,
1999, pp. 149-163.

\rm
\bibitem{Kon} Kontsevich, M. \bf Intersection Theory on the Moduli Space of Curves and the Matrix Airy Function,
\it functional analysis and it's applications, 1991,25:2, pp. 50-57.

\bibitem{Pen} Penner R.C. \bf The decorated Teichmuller Space of punctured
surfaces, \it Comm. Math. Phys., 113:2(1987), 299-340.

\bibitem{Chp} Shabat, G.B. \bf Combinatorial and Topological methods in the theory of algebraic curves,
\it Theses, Moscow State University, 1998.

\bibitem{K-CP} Shabat, G.B. \bf Complex analysis and dessins d'enfants,
\it in: "Complex Analysis in Modern Math.", "Phasis", Œoscow,  1998,
pp. 257-268.


\bibitem{P-CP} Shabat G.B., Voevodsky V. A. \bf Drawing curves over number
fields, \it The Grothendieck Festschrift, Birkhauser, 1990,  V.III.,
p.199-227.

\end{thebibliography}
\end{document}